\documentclass{article}

\usepackage{fixltx2e,fix-cm}

\usepackage{amssymb}
\usepackage{amsmath}
\usepackage{amsthm}[ntheorem]
\usepackage{graphicx,subcaption}
\usepackage{natbib}
\usepackage{makeidx}

\usepackage{thmtools}
\usepackage{stmaryrd,amsthm,url}
\usepackage{enumerate,enumitem}
\usepackage{etoolbox,mathtools}

\newtheoremstyle{theoremdd}
  {\topsep}
  {\topsep}
  {\upshape}
  {0pt}
  {\bfseries}
  {}
  { }
  {\thmname{#1}\thmnumber{ #2}\thmnote{ (#3)}}

\theoremstyle{theoremdd}

\newtheorem{prob}{Problem}

\declaretheorem[name=Example,qed={\lower-0.3ex\hbox{~\tiny$\blacksquare$}}]{ex}
\declaretheorem[name=Definition,qed={\lower-0.3ex\hbox{~\tiny$\blacksquare$}}]{defn}
\declaretheorem[name=Remark,qed={\lower-0.3ex\hbox{~\tiny$\blacksquare$}}]{remark}
\declaretheorem[name=Theorem,qed={\lower-0.3ex\hbox{~\tiny$\blacksquare$}}]{theorem}
\declaretheorem[name=Corollary,qed={\lower-0.3ex\hbox{~\tiny$\blacksquare$}}]{cor}
\declaretheorem[name=Lemma,qed={\lower-0.3ex\hbox{~\tiny$\blacksquare$}}]{lemma}

\numberwithin{remark}{section}
\numberwithin{theorem}{section}
\numberwithin{lemma}{section}
\numberwithin{equation}{section}
\numberwithin{defn}{section}
\numberwithin{figure}{section}
\numberwithin{table}{section}
\numberwithin{ex}{section}
\numberwithin{alg}{section}
\numberwithin{prob}{section}
\numberwithin{cor}{section}

\newcommand{\tcalf}{\tilde{\calf}}

\newcommand{\bals}{\begin{align*}}

\newcommand{\eals}{\end{align*}}

\newcommand{\absb}[1]{\ensuremath{\lvert#1\rvert}}
 
\newcommand{\absbb}[1]{\ensuremath{\left|#1\right|}}

\newcommand{\ben}{\begin{enumerate}}
\newcommand{\een}{\end{enumerate}}
\newcommand{\benr}{\begin{enumerate}[label=(\roman{*}),noitemsep,leftmargin=*]}
\newcommand{\benal}{\begin{enumerate}[label=(\alph{*}),noitemsep,leftmargin=*]}
\newcommand{\benar}{\begin{enumerate}[label=(\arabic{*}),noitemsep,leftmargin=*]}

\newcommand{\benrnls}{\begin{enumerate}[label=(\roman{*}),nolistsep,leftmargin=*]}
\newcommand{\benalnls}{\begin{enumerate}[label=(\alph{*}),nolistsep,leftmargin=*]}
\newcommand{\benarnls}{\begin{enumerate}[label=(\arabic{*}),nolistsep,leftmargin=*]}

\newcommand{\trace}{\mbox{trace}}

\newcommand{\bB}{\textbf{B}}

\newcommand{\bb}{\mathbf{b}}
\newcommand{\bU}{\mathbf{U}}

\newcommand{\bW}{\mathbf{W}}

\newcommand{\bz}{\mathbf{z}}

\newcommand{\diag}{\mbox{diag}}

\newcommand{\mthm}[1]{Theorem \ref{#1}}
\newcommand{\mlem}[1]{Lemma \ref{#1}}

\newcommand{\msec}[1]{Section \ref{#1}}
\newcommand{\mfig}[1]{Figure \ref{#1}}

\newcommand{\meq}[1]{Equation \eqref{#1}}

\newcommand{\bpr}{\begin{prob}}
\newcommand{\epr}{\end{prob}}

\newcommand{\bdm}{\begin{displaymath}}
\newcommand{\edm}{\end{displaymath}}

\newcommand{\beq}{\begin{equation}}
\newcommand{\eeq}{\end{equation}}

\newcommand{\bea}{\begin{eqnarray}}
\newcommand{\eea}{\end{eqnarray}}

\newcommand{\beas}{\begin{eqnarray*}}
\newcommand{\eeas}{\end{eqnarray*}}

\newcommand{\bdf}{\begin{defn}}
\newcommand{\edf}{\end{defn}}

\newcommand{\bex}{\begin{ex}}
\newcommand{\eex}{\end{ex}}

\newcommand{\bexa}{\begin{ex}}
\newcommand{\eexa}{\end{ex}}

\newcommand{\bexe}{\begin{exercise}}
\newcommand{\eexe}{\end{exercise}}

\newcommand{\bthm}{\begin{theorem}}
\newcommand{\ethm}{\end{theorem}}

\newcommand{\bmat}{\begin{bmatrix}}
\newcommand{\emat}{\end{bmatrix}}

\newcommand{\bproof}{\begin{proof}}
\newcommand{\eproof}{\end{proof}}

\newcommand{\blem}{\begin{lemma}}
\newcommand{\elem}{\end{lemma}}

\newcommand{\brem}{\begin{remark}}
\newcommand{\erem}{\end{remark}}

\newcommand{\bcor}{\begin{cor}}
\newcommand{\ecor}{\end{cor}}

\newcommand{\balg}{\begin{algorithm}}
\newcommand{\ealg}{\end{algorithm}}

\def\ng{%
  \setbox0=\hbox{-}%
  \vcenter{%
    \hrule width\wd0 height \the\fontdimen8\textfont3%
  }%
}

\newcommand{\bA}{\textbf{A}}
\newcommand{\bC}{\textbf{C}}

\newcommand{\bI}{\textbf{I}}

\newcommand{\bG}{\textbf{G}}

\newcommand{\bzero}{\mathbf{0}}

\newcommand{\calf}{{\cal F}}

\newcommand{\calx}{{\cal X}}

\newcommand{\caln}{{\cal N}}

\newcommand{\calt}{{\cal T}}

\newcommand{\real}{{\mathbb{R}}}

\newcommand{\integ}{{\mathbb{I}}}
\newcommand{\intege}{\bar{\integ}}

\newcommand{\eps}{{\epsilon}}

\newcommand{\bcb}{\begin{color}{blue}}
\newcommand{\bcr}{\begin{color}{red}}
\newcommand{\bcg}{\begin{color}{green}}
\newcommand{\ec}{\end{color}}

\newcommand{\pr}{{\prime}}

\newcommand{\sfrac}[2]{{#1}/{#2}}
\newcommand{\ssfrac}[2]{{#1}/\left({#2}\right)}
\newcommand{\sssfrac}[2]{\left({#1}/{#2}\right)}

\newcommand{\var}[1]{\mbox{var}[#1]}

\newcommand{\bsp}{\begin{sloppypar}}
\newcommand{\esp}{\end{sloppypar}}

\newcommand{\by}{\textbf{y}}
\newcommand{\bY}{\textbf{Y}}

\newcommand{\bx}{\textbf{x}}
\newcommand{\bX}{\textbf{X}}

\newcommand{\bn}{\textbf{n}}

\newcommand{\beps}{{\pmb{\epsilon}}}
\newcommand{\bbeta}{\pmb{\beta}}

\newcommand{\bu}{\pmb{u}}
\newcommand{\bv}{\pmb{v}}

\begin{document}

\title{A Stochastic Contraction Mapping Theorem}


\author{Anthony Almudevar \\
Department of Biostatistics and Computational Biology \\ University of Rochester}

\maketitle

\noindent\textbf{Abstract.}
In this paper we define contractive and nonexpansive properties for adapted stochastic processes $X_1, X_2, \ldots $ which can be used to deduce limiting properties. 
In general, nonexpansive processes possess finite limits while contractive processes converge to zero $a.e.$ Extensions  to multivariate processes are given. These properties may be used to model a number of important processes, including stochastic approximation and least-squares estimation of controlled linear models, with convergence properties derivable from a single theory.  The approach has the advantage of not in general requiring analytical regularity properties such as continuity and differentiability.   \\



\section{Introduction}

Let $\bX = \{ X_n : n \geq 0 \}$ be a sequence of random variables adapted to filtration $\tcalf = \{ \calf_n : n \geq 0 \}$. Following the Doob decomposition theorem we may always write 
\beq
X_n = X_{n-1} + A_n + \eps_n, \,\,\, n \geq 1,   \label{eq.intro.01}
\eeq
where $\eps_n = X_n - E[X_n \mid \calf_{n-1}]$ are $\calf_n$-measurable martingale differences, and  $A_n = E[X_n \mid \calf_{n-1}] - X_{n-1}$ defines a predictable process, in the sense that $A_n$ is $\calf_{n-1}$-measurable. We let $X_0$ be any suitable $\calf_0$-measurable initial value. We can think of $A_n$ as a control effector, and we will sometimes have $A_n = f_n(X_{n-1})$ for some sequence of deterministic mappings $f_n$, $n \geq 1$.   In this case, we can imagine setting   $\eps_n = 0$ in  \meq{eq.intro.01}, obtaining a deterministic iterative process. To fix ideas, suppose each $f_n$ has common fixed point $0 = f_n(0)$ and that there exists positive constants $\rho_n$ such that $0 \leq f_n(x)/x \leq \rho_n$ when $x \neq 0$.   In a manner similar to the Banach fixed point theorem, it can be shown that if $\prod_{n \geq 1} \rho_n = 0$, then $\lim_n X_n = 0$. We can refer to this as a contraction property, although we need not assume $f_n$ is formally a contraction mapping (or even Lipschitz continuous). If we restore $\eps_n$ to       \meq{eq.intro.01}, the question becomes how the additional stochastic variation affects the convergence properties. 

In \cite{almudevar2008approximate} (see also \cite{almudevar2014approximate}) a general theorem for such noisy iterative algorithms was developed, for which $X_n$ is defined on any Banach space $(\calx, \| \cdot \|)$ (not necessarily stochastic). For example, suppose operator $T_n$ possesses Lipschitz constant $\rho_n \leq \rho < 1$, and $x^*$ is the fixed point of each $T_n$. In this case the ``exact algorithm"   $X_n = T_n(X_{n-1})$, $n \geq 1$  converges to $x^*$. The ``approximate algorithm"  $X_n = T_n(X_{n-1}) + \eps_n$, $n \geq 1$, is then a noisy version of the exact algorithm.  Then if   $\limsup_n \| \eps_n \| = d$,  we may claim that  $\limsup_n \| X_n - x^*  \| = O(d)$, and if $d= 0$ then $X_n$ converges to $x^*$. In fact, the convergence rate of $\| X_n - x^*  \|$ is bounded by $O\left(  \max(\| \eps_n \|, (\rho+\delta)^n) \right)$ for any $\delta > 0$ \citep{almudevar2008approximate,almudevar2014approximate}.

The theory extends to nonexpansive operators $T_n$ ($\rho_n \leq 1$), especially ``weakly contractive" sequences which satisfy  $\prod_{n\geq 1} \rho_n = 0$, while, typically, the limit $\lim_n \rho_n = 1$ also holds.  In this case, stronger convergence results are obtainable when $(\calx, \| \cdot \|)$ is assumed to be additionally a Hilbert space (for example, adapted stochastic processes endowed with the $L_2$ norm). 

In this paper, we extend some of these ideas to filtered stochastic processes, with strong convergence in place of convergence in a Banach space. We will characterize precisely contractive or nonexpansive properties based on  the predictable process $A_n$ of \meq{eq.intro.01}. These force $X_n$ to behave like a supermartingale above 0, and a submartingale below. If $M_{1,n} = \sum_{i=1}^n \eps_i$ is a martingale of finite variation, it follows naturally that the process will possess a finite limit. However, if we then impose sufficient contractive properties, it can be shown that $X_n$ converges to 0 $a.e.$   


The main result is given in \msec{sec.main} for $X_n \in \real^1$. We introduce the nonexpansive and contractive properties, then \mthm{thm.main.nonexpansive} summarizes the important limiting properties.  

In \msec{sec.nonuniform} we consider the ``nonuniform contractive" case, by which we mean that the nonexpansive or contractive properties required of   \mthm{thm.main.nonexpansive} are now required to hold only outside a neighborhood of zero (Theorems \ref{thm.nonuniform.contraction.01} and \ref{thm.nonuniform.contraction.02}). 

In \msec{sec.multivariate},  \mthm{thm.main.nonexpansive} is extended to $\bX = \{ \bX_n \in \real^p : n \geq 0 \}$,   $p \geq 1$,  by replacing the constraint $0 \leq E[X_n \mid \calf_{n-1}]/X_{n-1} \leq \rho$ with  $\| E[\bX_n \mid \calf_{n-1}] \|/ \|\bX_{n-1}\| \leq \rho$,  $\| \bX_{n-1} \| \neq 0$. It is then quite straightforward to apply  \mthm{thm.main.nonexpansive} to $\| \bX_n \| \in \real^1$ (the conditions of which permit the process to be, for example, strictly positive).   
The argument is given in \mthm{thm.main.multivariate.01}. 

We include two important applications which can be modeled using the proposed theory. The Robbins-Monro stochastic approximation algorithm is a filtered stochastic process intended to converge to the solution to $g(\bx) = \bzero$ for a given mapping $g :\real^p \rightarrow \real^p$ \citep{robbins1951stochastic}. The univariate case will demonstrate the main theorem of \msec{sec.main}. Extension  to the nonuniform contractive model will significantly expand the class of $g$ to which stochastic approximation can be applied. The multivariate extension is then developed. This may be of some interest, since regularity conditions usually imposed on $g$ for $p \geq 1$ are  not needed by the methods introduced here    \citep{spall1992,pham2009contraction,lai2003stochastic,kushner2003stochastic,pham2009contraction}.

In \msec{sec.LS} we consider the controlled least-squares regression model introduced in \cite{christopeitHelmesAOS1980}.  We can directly derive the same conditions for strong convergence, while also proving that they are also necessary. In general, we can prove that the least-squares estimator always possess a finite limit, hence weak convergence implies strong convergence.  
 
Finally, we bring attention to a machine learning application proposed in \cite{hooker2022}, involving regularization of stochastic gradient boosted trees. This makes use of an earlier version of the methodology proposed here that appeared in an unpublished preprint written by this author. It demonstrates how the approach can be used to prove strong consistency where differentiability cannot be assumed.   

Throughout the manuscript proofs and a number of technical lemmas are given in the appendix.

\section{Main Result}\label{sec.main}

 Let $\bX=\{X_n : n \geq 0 \}$ be an $L_1$ stochastic process defined on probability space $(\Omega, \calf, P)$. Suppose there is a filtration $\tcalf$ defined defined by the $\sigma$-algebras $\calf_0  \subset \calf_1 \subset \calf_2 \subset \cdots \subset \calf$ to which $\bX$ is adapted. Define residuals
\begin{align}
\eps_n  =  X_n - E[X_n \mid \calf_{n-1}], \,\, n \geq 1, \label{eq.def.residuals.01}
\end{align}
and the partial sums
$
M_{s,t} = \sum_{i=s}^t \eps_i
$
for $t \geq s \geq 1$, with $M_{s,t} = 0$ whenever $t < s$. Then define the events
\begin{align*}
D^+_n  &= \{ X_n > 0 \},\,\,\, D^-_n  = \{ X_n < 0 \},\,\,\, D^0_n =  \{ X_n = 0 \}, \,\, n \geq 0, \\
D_n & =  D^+_n (D^+_{n-1})^c \cup D^-_n (D^-_{n-1})^c \cup D^0_n(D^0_{n-1})^c, \,\, n \geq 1, 
\end{align*}
and the compound events $D^+_{s,t}  =  \bigcap_{i=s}^t D^+_i$, $D^+_{s,\infty}  =  \bigcap_{i \geq s} D^+_i$, $0 \leq s \leq t$.
Define similarly $D^-_{s,t}$, $D^-_{s,\infty}$, $D^0_{s,t}$ and $D^0_{s,\infty}$. We may then construct the increasing sequence of stopping times $T_1,T_2,\ldots$
\begin{displaymath}
T_j \leq t \Leftrightarrow  \sum_{i=1}^t I\{D_i\} \geq j
\end{displaymath}
for $j \geq 1$. Then $T_1,T_2,\ldots$ may be interpreted as the times at which sequence $\bX$ changes sign, changes from zero to a nonzero value or changes from a nonzero value to zero. We can refer to any such event as a crossing. Note that $T_j$ may be infinite, so it will be useful to define
$
N_T = \sup\{ j \geq 1 : T_j < \infty\}.
$
Then let
\begin{displaymath}
W_j = \left\{ \begin{array}{ccc}  \sup_{T_j \leq i < T_{j+1}}
\absb{X_i} & ; & T_j < \infty \\
0 & ; & T_j = \infty \end{array} \right.
\end{displaymath}
for $j \geq 1$. We will make use throughout the paper of the following sequence:
\beq
U_n = E[X_n  \mid  \calf_{n-1}] I\{X_{n-1} = 0\},\,\,\, n \geq 1. \label{eq.def.u.01}
\eeq
Our model will be based on the following assumptions:
\begin{itemize}[noitemsep]
\item[(A1)] There exist nonnegative constants $\alpha_n \geq 0$, $n \geq 1$, such that 
$$
 0 \leq \frac{E[X_n  \mid  \calf_{n-1}]}{X_{n-1}} I\{X_{n-1} \neq 0\}  \leq 1+\alpha_n, \,\,\, a.e.
$$
and $\sum_{i \geq 1} \alpha_i < \infty$. 
\item[(A2)] There exist nonnegative constants $0 \leq k_n \leq 1$, $n
\geq 1$, such that
$$
 0 \leq  \left(\sfrac{E[X_n  \mid  \calf_{n-1}]}{X_{n-1}}\right) I\{X_{n-1} \neq 0\}  \leq k_n, \,\,\, a.e.
$$
and $\sum_{i \geq 1} (1-k_i) = \infty$ (equivalently, $\prod_{i \geq s} k_i = 0$ for all $s \geq 1$).  
\item[(A3)]$\lim_{n\rightarrow\infty} U_n = 0 \,\, a.e.$
\item[(A4)] The partial sums $M_{1,n}$ possess a finite limit $a.e.$ as $n\rightarrow\infty$.
\end{itemize}

Assumption (A1) is analogous to a nonexpansive condition on $E[X_n  \mid  \calf_{n-1}]$, and is related to the ``almost martingale"  
introduced in \cite{robbins1971convergence}.  Assumption (A2) is the stronger, or contractive, version of (A1). The central result of this section is that under assumptions (A3) and (A4) the process defined by (A1) converges to a finite limit, while that defined by (A2) converges to 0.

In some applications, (A3) will be a natural extension of (A1) or (A2) to the case $X_{n-1} = 0$. That is, we may be able to  claim  $0 \leq E[X_n  \mid  \calf_{n-1}] \leq X_{n-1}(1+\alpha_n)$ for $X_{n-1} \geq 0$ and $0 \geq E[X_n  \mid  \calf_{n-1}] \geq X_{n-1}(1+\alpha_n)$ for $X_{n-1} \leq 0$. In this case (A3) holds trivially. However, we find that our theory can be extended to a broader class of interesting models by permitting somewhat more flexibility when conditioning on $\{ X_{n-1} = 0 \}$. Finally, since $M_{1,n}$ $n \geq 1$ is a martingale, conditions under which (A4) holds are well known. For example, if $M_{1,n}$ is an $L_2$ martingale, by Theorem 4.5.2 of \cite{durrett2019probability}, we may claim 
\beq
\sum_{i=1}^\infty E[\eps_i^2 \mid \calf_{i-1} ] < \infty \,\,\, a.e. \mbox{   implies   } \lim_n M_{1,n}  \mbox{  exists and is finite $a.e.$} \label{eq.durret.4.5.2}
\eeq


 
 
We are now give the main results. We first consider the nonexpansive model. Given (A1), (A3) and (A4) we may claim that $\bX$ possesses a finite limit. 
When (A2) holds as well we have convergence to zero. It will be useful, however, to consider the case for which (A3) need not hold. The results are essentially partitioned into the cases $\{ N_T < \infty\}$ and $\{ N_T = \infty\}$, since the analysis differs between these cases in some important ways. For example, on $\{ N_T = \infty\}$ the nonexpansive property (A1) is sufficient for convergence to zero.


\begin{theorem}\label{thm.main.nonexpansive} 
Suppose  $\bX$ satisfies (A1) and (A4). Then the following statements hold.
\benr
\item
$X_n I\{ N_T < \infty\}$ possesses a finite limit $a.e.$ as $n\rightarrow\infty$. 
\item
If in addition (A2) holds then $X_n I\{ N_T < \infty\}$ converges to zero $a.e.$
\item
The inequality $\limsup_n \absb{X_n} I\{ N_T = \infty \} \leq  \limsup_n \absb{U_n}$ holds $a.e.$
\item
If in addition (A3) holds, then $\lim_n X_n I\{ N_T = \infty \} = 0$ $a.e.$ \qedhere
\een
\end{theorem}

 




\subsection{The Strong Law of Large Numbers as a Contractive Process}\label{sec.slln} 

We now discuss the relationship between \mthm{thm.main.nonexpansive} and the strong
law of large numbers. We suppose that there is an $L_1$ sequence of random
variables $\bY=\{Y_n : n \geq 0\}$ on a probability space $(\Omega, \calf, P)$ such that $S_n = \sum_{i=0}^n Y_i$, $n \geq 0$, is a martingale with respect to some filtration $\tcalf$.

Suppose $a_0, a_1, a_2, \ldots$ is a nondecreasing sequence of positive real numbers such that $a_n \rightarrow_n \infty$. Then for any $s \geq 1$, 
\beq
\lim_{n\rightarrow\infty} \prod_{i =s}^n \frac{a_{i-1}}{a_i} = \lim_{n\rightarrow\infty} \frac{a_{s-1}}{a_n} = 0. \label{eq.kronecker.01}
\eeq
We are interested in conditions under which
$
\sfrac{S_n}{a_n} \rightarrow_n  0 
$
$a.e.$ We first note that, following  \meq{eq.kronecker.01}, the sequence $X_n = S_n/a_n$ satisfies assumption (A2), since
$$
\frac{E[X_n  \mid  \calf_{n-1}]}{X_{n-1}} I\{X_{n-1} \neq 0\}  = 
\frac{S_{n-1}/a_n}{S_{n-1}/a_{n-1}} I\{X_{n-1} \neq 0\} 
 \leq \frac{a_{n-1}}{a_n},\,\,\, n \geq 1. 
$$
Then (A3) holds since
$
E[X_n  \mid  \calf_{n-1}]I\{X_{n-1} = 0\}  =  (\sfrac{S_{n-1}}{a_n}) I\{X_{n-1} = 0\} =  0.
$
The residual becomes
$
\eps_n  =  X_n - E[X_n \mid \calf_{n-1}] 
     =  \sfrac{S_n}{a_n} - \sfrac{S_{n-1}}{a_n} 
     =  \sfrac{Y_n}{a_n}.
$
Direct application of \mthm{thm.main.nonexpansive}  gives
\begin{theorem}
Let $S_n= \sum_{i=0}^n Y_i$ be an $L_1$ martingale with respect to filtration $\tcalf$. If  $a_n\uparrow_n\infty$, and if $\sum_{i=1}^\infty Y_i/a_i$ is a convergent series $a.e.$ then $S_n/a_n \rightarrow_n 0$ $a.e.$
\end{theorem}
\noindent But this is simply a restatement of Kronecker's lemma, which states that if $a_n \uparrow_n \infty$ then for any sequence of real numbers $x_1, x_2, \ldots$ the convergence of the series $\sum_{i=1}^\infty x_i/a_i$ implies $\lim_n (1/a_n) \sum_{i=1}^n x_i = 0$. 


\subsection{Stochastic Approximation}\label{sec.sa.01}

The stochastic approximation algorithm  is a method used to determine the solution $x^* \in E_g \subset \real$ to
\beq
g(x) = 0 \label{eq.sa.gx.01}
\eeq
where mapping $g:E_g \rightarrow\real$ can only be evaluated with noise. In \cite{robbins1951stochastic} the iterative algorithm 
\begin{equation}
X_n = X_{n-1} - \alpha_nU_n, \,\,\, n \geq 1 \label{Algorithm}
\end{equation}
was proposed, giving conditions under which $X_n$ converges in $L_2$  to $x^*$. Here, $\alpha_n$ is a sequence of positive constants, and the sequence $X_n$ is adapted to filtration  $\tcalf$ and constructed so that $E[U_n \mid \calf_{n-1}] = g(X_{n-1})$. 

Strong convergence has been since been established under a variety of conditions \citep{blum1954approximation, robbins1971convergence, ljung1978strong,   lai2003stochastic, kushner2003stochastic, benveniste2012adaptive}. Our first task with respect to stochastic approximation will be to show how it may be interpreted as a stochastic contraction process. 

Given a sequence $\bX$ constructed by the iterations of  \meq{Algorithm} we define the following
assumptions.
\begin{itemize}[noitemsep]
\item[(B1)] The mapping $g:E_g  \rightarrow\real$,  possesses root $x^* \in E_g$.
\item[(B2)] $P(X_n \in E_g) = 1$ for all $n \geq 1$.
\item[(B3)] There exists $0 < m \leq M < \infty$ such that
$
m \leq \ssfrac{g(x)}{x - x^*} \leq M
$
for all $x \in E_g - \{x^*\}$.
\item[(B4)] $\alpha_n \geq 0$,  $\lim_{n \rightarrow\infty} \alpha_n = 0$.
\item[(B5)] The series $\sum_{i=0}^\infty \alpha_i \left\{ U_i - E[U_i \mid \calf_{i-1}] \right\}$ converges $a.e.$
\item[(B6)] $\sum_{i=1}^\infty \alpha_i = \infty$.
\item[(B7)] $E[U_n \mid \calf_{n-1}] = g(X_{n-1}) $ $a.e.$ for $n\geq 1$.
\end{itemize}

Assumption (B1) does not require that $x^*$ be in the interior of $E_g$. Nor does  it state that the root is unique, but this is forced by (B3).  Then, 
if finite nonzero left and right derivatives of $g$ exist at $x^*$, (B3) will hold on some open neighborhood of $x^*$, if not the entire domain of $g$. 

It should be noted that, following \meq{eq.durret.4.5.2},    if the variance of $U_n$ conditional on $\calf_{n-1}$ is bounded, then  (B5) may be replaced by $\sum_{i \geq 1} \alpha^2_i < \infty$, which, with (B6), are the conditions on the sequence $\alpha_n$ commonly given in the literature

We now state the main theorem of this section. 

\begin{theorem}\label{thm.rm.01}
If (B1)-(B7) hold then the process $\bX$ defined by \meq{Algorithm} converges to $x^*$ $a.e.$
\end{theorem}

It is interesting to note that the stochastic approximation algorithm will converge if $g(x)$ can be perfectly evaluated (in contrast with simulated annealing, which relies on stochastic variation for convergence).  In this case \meq{Algorithm} becomes
$$
x_n = x_{n-1} - \alpha_n g(x_{n-1}) = T_n(x_{n-1}), \,\,\, n \geq 1,
$$
where $T_n$ is a sequence of operators with Lipschitz constant $\rho_n  = (1 - m\alpha_n)$ for all large enough $n$.  Given (B6)  we have  $\prod_{i \geq 1} \rho_i  = 0$,  so that although the Lipschitz constants are not bounded away from one, they approach one from below slowly enough to allow contraction to force convergence. Two conditions on $\alpha_n$ are commonly associated with stochastic approximation, in particular, (a) $\sum_{i=1}^\infty \alpha_i = \infty$; and (b) $\sum_{i=1}^\infty \alpha_i^2 <  \infty$. As can be seen, condition (a) forces sufficient contraction to ensure convergence for the noiseless algorithm, while condition (b) ensures that the cumulative effect of noise vanishes in the limit (as a consequence of contraction property). For more on this point of view, see \cite{almudevar2008approximate,almudevar2014approximate}.

\section{Nonuniform Contractions}\label{sec.nonuniform}

\begin{figure}
\centering
\begin{subfigure}{0.465\textwidth}
\includegraphics[width=\textwidth, viewport = 0 15 475 450, clip]{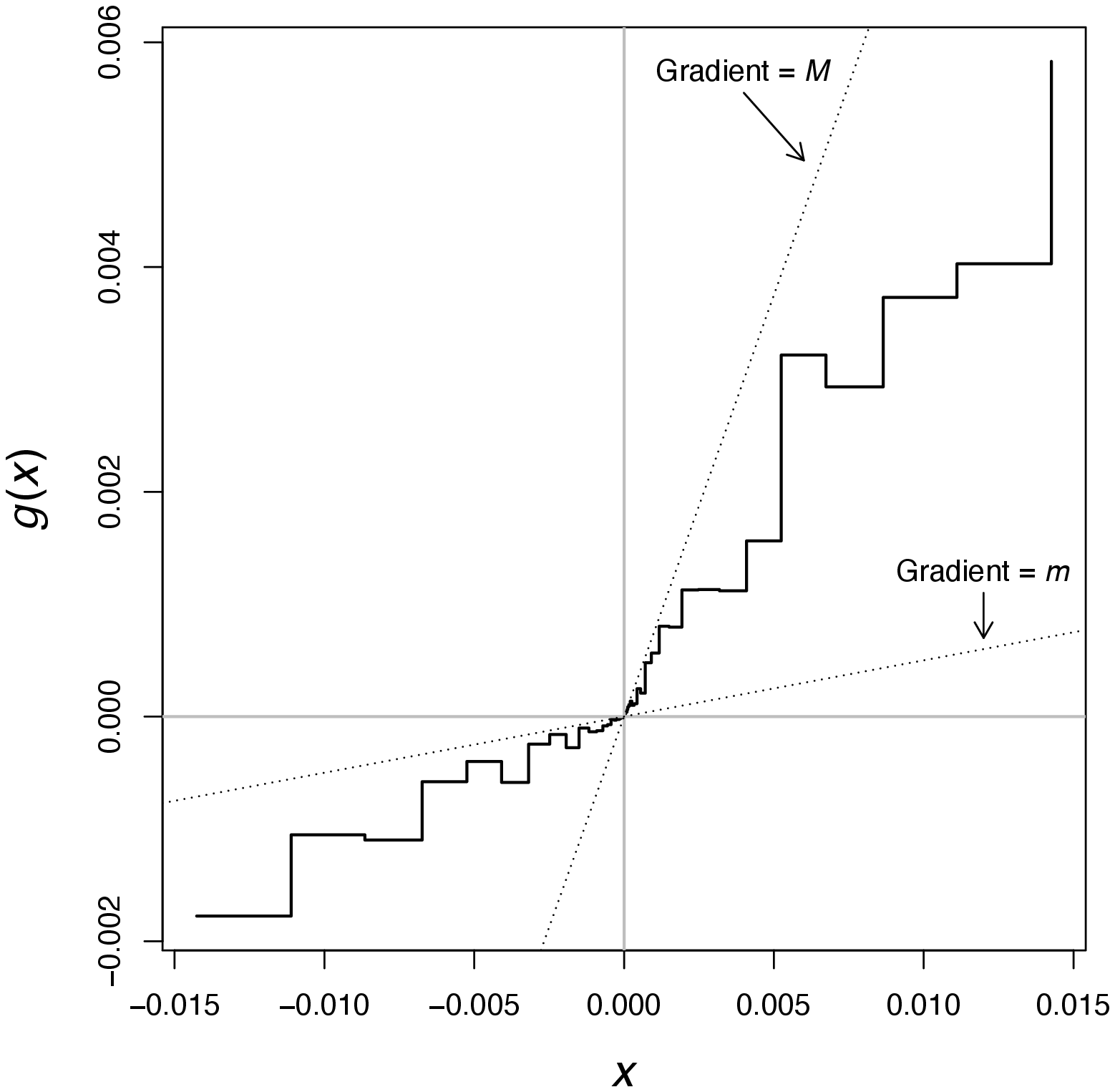}
\caption{}\label{scfig1}
\end{subfigure}
\hfill
\begin{subfigure}{0.475\textwidth}
\includegraphics[width=\textwidth, viewport = 0 15 475 450, clip]{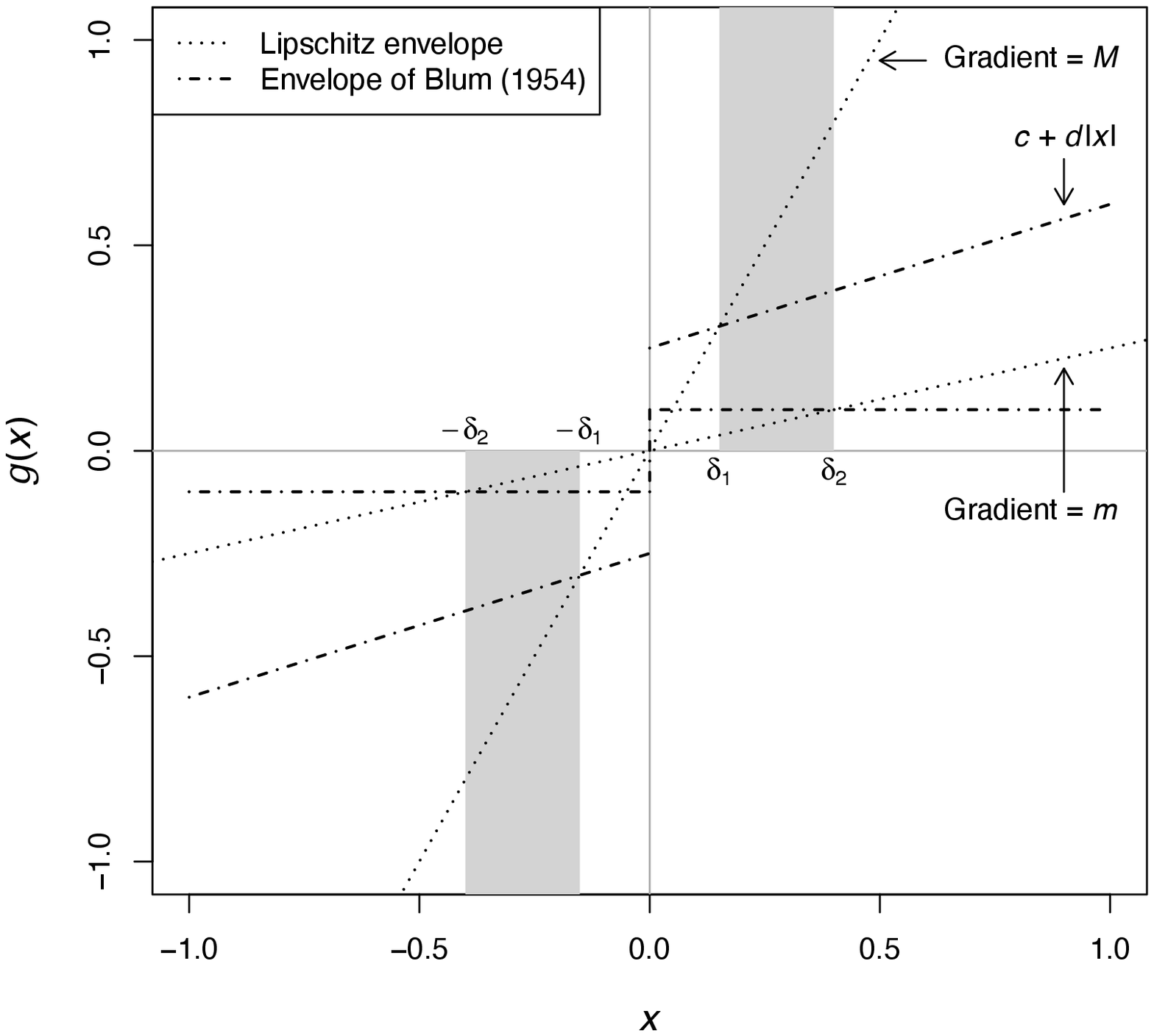}
\caption{}\label{scfig2}
\end{subfigure}
\caption{Figure (a) shows the Lipschitz envelope defined by gradients $M$ and $m$ of condition (B3), which  essentially represent sufficient conditions on the function $g(x)$ of Equation \eqref{eq.sa.gx.01} for convergence of the stochastic approximation algorithm using \mthm{thm.main.nonexpansive}. Note that monotonicity and continuity are not required. Figure (b) gives a comparison of the regularity conditions for objective function $g(x)$ based on  \mthm{thm.main.nonexpansive}, and the nonuniform contraction model of \mthm{thm.nonuniform.contraction.01}. The Lipschitz envelope shown in \mfig{scfig1} is superimposed here, along with a graphical representation of the comparable envelope forced by regularity conditions (C1)-(C3)  introduced in \cite{blum1954approximation}. The important feature is that for any $\delta >0$ we may construct an upper Lipschitz envelope which bounds above the upper envelope of (C1) for all $\absb{x} > \delta$, and for any $0 < \delta_1 < \delta_2 < \infty$, we may construct a lower Lipschitz envelope which bounds  below the lower envelope implied by (C3) (see the gray shaded area).}
\end{figure}

\mfig{scfig1} demonstrates that the conditions imposed on $g(x)$ by \mthm{thm.rm.01} are quite general. Essentially, $g(x)$ must be contained within a Lipschitz envelope, but otherwise need not be monotone or continuous  (for convenience we will assume $x^* = 0$).

However, a variety of conditions on $g$ have been derived in the literature. In the seminal paper \cite{robbins1951stochastic}, alternative conditions for $L_2$ convergence   are given, in particular, $g(x) \leq -\delta$ for $x < 0$ and  $g(x) \geq \delta$ for $x > 0$, for some $\delta > 0$ (Theorem 1); and $g(0) = 0$, $g(x)$ is nondecreasing and $g^\prime(0) > 0$ (Theorem 2).  In \cite{blum1954approximation}, the following regularity conditions are offered:
\begin{itemize}[noitemsep]
\item[(C1)] $\absb{g(x)} \leq c + d\absb{x}$ for some constants $c,d \geq 0$; 
\item[(C2)] $g(x) < 0$ for $x < 0$ and $g(x) > 0$ for $x > 0$;
\item[(C3)] $\inf_{\delta_1 \leq \absb{x} \leq \delta_2} \absb{g(x)} > 0$ for every pair of numbers $\delta_1, \delta_2$ for which $0 < \delta_1 < \delta_2 < \infty$. 
\end{itemize}
While conditions (C1)-(C3) are considerably weaker than the envelope bound of  Figure \ref{scfig1}, they can be usefully compared. If (C1) holds, then some Lipschitz upper bound (Figure \ref{scfig1}) holds for $\absb{x} > \delta$. Then, with (C2), we can show that (A1) will hold \emph{outside} a neighborhood $\caln$ of the origin. Then if (C3) holds we may show that (A2) holds for  $\delta_1 \leq \absb{x} \leq \delta_2$, for any positive constants $\delta, \delta_1 < \delta_2$. In this case, it will be possible to show, first, that $\bX$ will be bounded in the limit. In this case, (A2) need only hold for a bounded subset of $E_g$, and we will be able to show that $\bX$ approaches $\caln$ $a.e.$, essentially ``fattening" the origin. If the neighborhood $\caln$ can then be made arbitrarily small, we obtain strong convergence. 

Thus, we can say that, in the context of stochastic approximation, (C1)-(C3) define models that are contractive outside any neighborhood of the origin, hence it seems reasonable to refer to such a model as ``nonuniformly contractive". Below, we formally define such a model.  We again define an adapted sequence $\bX$, on which we may, as an alternative to (A1)-(A4), impose the following assumptions.   
\begin{itemize}[noitemsep]
\item[(D1)] For $\delta > 0$ (A1) holds for $\absb{ X_{n-1} } > \delta$ $a.e.$
\item[(D2)] For $0 < \delta_1 < \delta_2 < \infty$ (A2) holds for $\absb{ X_{n-1} } \in (\delta_1, \delta_2)$ $a.e.$ 
\item[(D3)]  $\absb{X_{n-1}} \in (0, \delta]$ implies  $\absb{ E[X_n  \mid \calf_{n-1}] } \leq \kappa < \infty$ $a.e.$ 
\item[(D4)] For all $\delta > 0$, (D3) holds for constants $\kappa_\delta$ which  satisfy the limit $\lim_{\delta\downarrow 0} \kappa_\delta = 0$. 
\end{itemize}
The key to the argument is to consider the process derived from $\bX$:
\begin{equation}
X^{\delta+\tau}_n = X_n I\{ \absb{ E[ X_n \mid \calf_{n-1}] } \geq \delta + \tau \}, \,\,\, n \geq 1. \label{eq.def.x.delta.01}    
\end{equation}
for given $\delta, \tau > 0$ (it will prove to be more intuitive to introduce $\delta$ and $\tau$ separately in our notation). We denote this process $\bX^{\delta+\tau}$. We may define the  residuals in the same manner as in \meq{eq.def.residuals.01}:
\begin{equation}
\eps_{n}^{\delta+\tau} = X_{n}^{\delta+\tau} - E[X_{n}^{\delta+\tau} \mid \calf_{n-1}], \,\,\, n \geq 1, \label{eq.def.x.delta.residuals.01}    
\end{equation}
Finally, following \meq{eq.def.u.01} we may condition the process on $\{X^{\delta+\tau}_{n-1} = 0 \}$:
\begin{equation}
U^{\delta+\tau}_{n}  = E[X_{n}^{\delta+\tau} \mid \calf_{n-1}] I \{  X^{\delta+\tau}_{n-1} = 0 \}. \label{eq.u.derived.process.01}
\end{equation}

We briefly outline the remaining strategy.  Suppose we can assume that $\eps_n \rightarrow_n 0$ $a.e.$, where $\eps_n$  are the residuals of the original process $\bX$. This will be the case under (A4).  Then under (D1), there exists $wp1$ a finite integer $N_0$ such that $\bX^{\delta+\tau}$ satisfies (A1) for all $n \geq N_0$. Accordingly, we define the offset process $\bX^{\delta+\tau}_{N_0} = \{ X^{\delta+\tau}_{N_0 + n} : n \geq 0\}$, so we may claim that $\bX^{\delta+\tau}_{N_0}$ satisfies (A1).  The remaining argument is of the same style, that is, by assuming the various conditions (D1)-(D4) hold for $\bX$, we may verify that the corresponding properties (A1)-(A4) hold for    $\bX^{\delta+\tau}_{N_0}$. Then the convergence properties of $\bX$ follow from those of $\bX^{\delta+\tau}_{N_0}$.  The first theorem follows.

\begin{theorem}\label{thm.nonuniform.contraction.01}
Suppose $\eps_n \rightarrow_n 0 $ $a.e.$  Suppose (D1) holds for some $\delta > 0$. Consider the derived process of Equation \eqref{eq.def.x.delta.01}. Then 
\benr
\item
For all small enough $\tau > 0$, $wp1$ there exists finite $N_0$  such that $\bX^{\delta+\tau}_{N_0}$ satisfies (A1) for all $n \geq 0$.
\item
If in addition, $\bX$ satisfies (D2) for $\delta_1 \geq \delta$, then $\bX^{\delta+\tau}_{N_0}$  satisfies (A2) whenever $\absb{X^{\delta+\tau}_{N_0 + n-1}} < \delta_2$. 
\item
If in addition (D3) holds, then $wp1$ 
\begin{equation}
 U^{\delta+\tau}_{N_0+n} \leq (\absb{U_{N_0+n}} + \delta  + 2\tau + \kappa), \,\,\, n \geq 0. \label{eq.thm6.(iii)}
\end{equation}   
\item
If (A4) holds for $\bX$, then it also holds for $\bX_{N_0}^{\delta+\tau}$. \qedhere
\een
\end{theorem}

\mthm{thm.nonuniform.contraction.01} allows us to deduce some convergence properties of $\bX^{\delta+\tau}$. If (D1) and (A4) hold for $\bX$ then (A1) and (A4) hold for  $\bX^{\delta+\tau}_{N_0}$. As a consequence, by \mthm{thm.main.nonexpansive} either   $X_{n}^{\delta+\tau}$ converges to a finite limit, or 
$
\limsup_{n \rightarrow \infty} \absb{X_{n}^{\delta+\tau}} \leq \limsup_{n \rightarrow \infty} \absb{U^{\delta+\tau}_{n}}.
$
The final step is to deduce the implications of this for the limit of $\bX$ itself.

\begin{theorem}\label{thm.nonuniform.contraction.02} Suppose $\bX$ satisfies  (D1) for all $\delta > 0$ $a.e.$, and in addition (D4), (A3) and (A4).  Then $\bX$ possesses a finite limit $a.e.$ If in addition $\bX$ satisfies (D2) for all pairs  $0 < \delta_1 < \delta_2 < \infty$ then $\bX$ converges to zero $a.e.$
\end{theorem}

\subsection{Stochastic Approximation and Nonuniform Contraction}\label{sec.blum.sa}

We continue our discussion of stochastic approximation. Consider assumptions (B1)-(B7) (as before, assuming $x^* = 0$). Then replace assumption (B3), which is placed on $g(x)$, with (C1)-(C3). 
Fix any $\delta > 0$. By (C1) $\absb{g(x)} \leq c + d \absb{x}$. With (C2) we may conclude
$$
0 \leq \sfrac{g(x)}{x} \leq \sssfrac{c}{\delta} + d, \,\,\, \mbox{for all}\,\,\, \absb{x} > \delta.
$$
Since $\alpha_n \rightarrow_n 0$ we have, for all large enough $n$, following \meq{eq.thm.rm.01.01}, 
$$
0 \leq \left(1 - \alpha_n \left( \frac{c}{\delta} + d\right) \right)  I\{X_{n-1} \neq 0\}  \leq \frac{ E[X_n  \mid  \calf_{n-1}] }{X_{n-1}} I\{X_{n-1} \neq 0\} \leq 1, 
$$
so that assumption (D1) will be satisfied. Thus, (C1)-(C2) force the nonexpansive property.  

To verify (D4), note that by (C1)-(C2) $\absb{E[X_n  \mid  \calf_{n-1}]} \leq \absb{X_{n-1}}$ for $\absb{X_{n-1}} \leq \delta$ for all large enough $n$. Thus, (D4) holds by setting $\kappa_\delta = \delta$.   

We next show that (C3) forces the contractive property. Fix $0 < \delta_1 < \delta_2 < \infty$. Define the quantity
$$
K_{\delta_1, \delta_2} = \inf_{\delta_1 \leq \absb{x} \leq \delta_2} g(x)/x.
$$
Given (C2), under (C3) $K_{\delta_1, \delta_2} > 0$, so that, for all large enough $n$, 
$$
0 \leq \frac{ E[X_n  \mid  \calf_{n-1}] }{X_{n-1}} I\{X_{n-1} \neq 0\} \leq \left(1 - \alpha_n K_{\delta_1, \delta_2} \right).
$$
Following the proof of \mthm{thm.rm.01}, we may then conclude that (D2) holds for all finite pairs  $\delta_1, \delta_2$. Thus, by \mthm{thm.nonuniform.contraction.02},
the stochastic approximation algorithm will converge when (B3) is replaced by (C1)-(C3).

\section{Multivariate Processes}\label{sec.multivariate}

We now extend \mthm{thm.main.nonexpansive} to processes in $\real^p$, replacing the ratio in (A1)  and (A2) with the absolute value of the ratio. Suppose $\bX=\{\bX_n \in\real^p ; n\geq 0\}$ is stochastic process adapted to filtration $\tcalf$.   The residual vector is $\beps_n = \bX_n - E[\bX_n  \mid  \calf_{n-1}]$, $n \geq 1$. Then $E[\beps_n \mid \calf_{n-1}] = 0$, so that as before $\beps_n$, $n \geq 1$ form vectors of martingale differences. Denote the components of  $\bX_n$ and $\beps_t$ by   $\bX_n(t) \in\real$ and $\beps_n(t)\in\real$, $t \in \calt = \{1, \ldots, p\}$.  Consider the following assumptions, which may be taken as
multivariate extensions of (A1)-(A4) (here, $\| \cdot \|$ is the Euclidean norm).

\begin{itemize}[noitemsep]
\item[(A5)] There exist nonnegative constants $\alpha_n \geq 0$, $n \geq 1$, such that 
\begin{displaymath}
 \frac{\| E[\bX_n  \mid  \calf_{n-1}] \|}{\| \bX_{n-1}\|} I\{ \|\bX_{n-1} \| \neq 0\}
 \leq 1+\alpha_n \,\, a.e.
\end{displaymath}
and $\sum_{i \geq 1} \alpha_i < \infty$. 
\item[(A6)] There exists positive constants $k_n \leq 1$, $n
\geq 1$, such that
\begin{displaymath}
 \frac{ \| E[\bX_n  \mid  \calf_{n-1}] \|}{ \|\bX_{n-1} \| } I\{ \|\bX_{n-1} \| \neq 0\} \leq
 k_n \,\, a.e. 
\end{displaymath}
and
$\sum_{i \geq 1} (1-k_i) = \infty$ (equivalently, $\prod_{i \geq s} k_i = 0$ for all $s \geq 1$).  
\item[(A7)]$\lim_{n\rightarrow\infty} E[\bX_n \mid \calf_{n-1}] I\{\bX_{n-1} = 0\} = 0
\,\, a.e.$
\item[(A8)] There exists a sequence of constants $v_n \geq 0$, $n \geq 1$ such that $\sum_{i=1}^\infty v_i < \infty$  and
$\max_{t \in\calt} \var{\beps_n(t) \mid \calf_{n-1}} \leq v_n$ for all $n \geq 1$. 
\end{itemize}
Assumptions (A5) and (A6) are strict generalizations of (A1) and (A2) even for $p=1$, since the contraction implied by $E[\bX_n \mid \calf_{n-1}]$ is no longer required
to be of the same sign as $\bX_n$. Assumptions (A7)-(A8) are directly comparable to (A3)-(A4), respectively.

The strategy will be to apply  \mthm{thm.main.nonexpansive} directly to $\| \bX_n \|$, $n \geq 0$. 
\begin{theorem}\label{thm.main.multivariate.01}
Suppose $\bX=\{\bX_n; n \geq 0\}$ is an $L_2$ process in $\real^p$, $p
\geq 1$ which is adapted to filtration $\tilde{\calf}$.
\benr
\item If (A5), (A7), (A8) hold then $\lim_{n\rightarrow\infty} \| \bX_n \|$
exists  and is finite $a.e.$
\item If (A6), (A7), (A8) hold then $\lim_{n\rightarrow\infty} \| \bX_n\| = 0$ $a.e.$ \qedhere
\een
\end{theorem}

\subsection{Multivariate Stochastic Approximation}

In this section we extend the stochastic approximation algorithm to the multivariate case. In particular, we consider the problem of determining the solution to
 $\bx^* \in E_g \subset \real^p $ to
$
g(\bx) = \bzero
$
where $g:E_g \rightarrow\real^p$ is only observable with noise. The extension to $\bx^* \in \real^p$ has been widely considered in the literature \citep{ljung1978strong,   spall1992,lai2003stochastic,kushner2003stochastic}. Regularity conditions usually impose differentiability assumptions on $g(\bx)$, although in 
\msec{sec.sa.01} and \msec{sec.blum.sa} it can be seen that such conditions play little role for the univariate case. We show in this section how \mthm{thm.main.multivariate.01} can be used to extend this nonsmooth case to $\real^p$ (we do not consider here the nonuniform contractive model, although presumably it would also apply to this application). 

We keep more or less intact conditions (B1)-(B7), with (B3) becoming a type of positive definiteness assumption on $g(\bx)$. Assuming $\bx^* = \bzero$, the algorithm retains form 
\begin{equation}
\bX_n = \bX_{n-1} - \alpha_n \bU_n,\,\,\, n \geq 1, \label{Algorithm2}
\end{equation}
here interpreted as a multivariate process, where $\bX_n, \bU_n  \in \real^p$, and $\alpha_i$ is a positive constant. Again, $\bX_n$, $n \geq 0$,  is adapted to some filtration $\tcalf$. We replace (B1)-(B7) with the following assumptions (some are identical to the corresponding assumption in  (B1)-(B7), but are included for convenience).  
\begin{itemize}[noitemsep]
\item[(B1a)] The mapping $g:E_g \rightarrow\real^p$, $E_g \subset \real^p$, possesses root $\bzero \in E_g$.
\item[(B2a)] $P(\bX_n \in E_g) = 1$ for all $n \geq 1$.
\item[(B3a)] There exists $0 < m \leq M < \infty$ such that $g$ satisfies 
$$
m \leq (g(\bx) \circ \bx)/\|\bx\|^2\,\,\,\mbox{and}\,\,\, \|g(\bx)\|/\|\bx\| \leq M 
$$ 
for all $\bx \in E_g - \{\bzero \}$, where $\bx \circ \bx^\pr$ is the inner product on $\real^p$. 
\item[(B4a)]  $\alpha_n \geq 0$, $\lim_n \alpha_n = 0$. 
\item[(B5a)] $\sum_{i=1}^\infty \alpha_i^2 \var{\bU_i(t) \mid \calf_{i-1}} < \infty$ $a.e.$ for all $t \in \calt$.
\item[(B6a)] $\sum_{i=1}^\infty \alpha_i = \infty$.
\item[(B7a)] $E[\bU_n \mid \calf_{n-1}] = g(\bX_{n-1}) $ $a.e.$ for $n\geq 1$.
\end{itemize}

Under these conditions the multivariate stochastic approximation algorithm remains strongly convergent. 

\begin{theorem}\label{thm.rm.02}
If (B1a)-(B7a) hold then the process $\bX_n$, $n\geq 1$, defined by \meq{Algorithm2} converges to $\bzero$ $a.e.$
\end{theorem}

\section{Least Squares Regression and Control}\label{sec.LS}

We consider the example of multivariate linear regression, possibly controlled. We are given probability measure space $(\Omega, \calf, P)$ on which $\tcalf$ is a filtration.  Suppose $\{\bx_i = (x_{i1}, \ldots, x_{ip}) ; i\geq 1\}$ is a sequence of random $p\times 1$ vectors, where $\bx_i$ is $\calf_{i-1}$-measurable.  Let $\bX_n$ be the $n\times p$ matrix with $i$th row equal to $\bx^T_i$. Suppose $u_1, u_2,\ldots$ is a sequence of zero mean random variables. Define $n\times 1$ vector $\bu_n = [u_1  \ldots  u_n ]^T$. We assume $\bu_n$  is $\calf_{n}$-measurable,  and defines a sequence of martingale differences, that is, $E[u_n \mid \calf_{n-1}] = 0$, $n \geq 1$.  Let $\bbeta = [\beta_1 \ldots \beta_p]^T$ be a fixed $p\times 1$ vector. Define $d^2_{n,t} = \sum_{i=1}^n x_{it}^2$, $d^2_{0,t} = 0$ and $\bA_n = \bX^T_n\bX_n$.  Then consider the nested sequence of linear models
\begin{equation}
\bY_n = \bX_n \bbeta + \bu_n, \,\, n \geq 1. \label{linearmodel}
\end{equation}
If $\bX_n$ contains $p$ linearly independent rows, $\bX_{n^\prime}$ does as well for any $n^\prime > n$. It will therefore be reasonable to assume that $wp1$ there exists finite $n_0$ such that $\bA_n$ is nonsingular for $n \geq n_0$. 

The ordinary least squares estimate of $\bbeta$ based on the $n$th model is well known to be $\bb_n = \bA_n^{-1}\bX^T_n\bY_n$. We then center the estimates by setting $\bar{\bb}_n = \bb_n - \bbeta = \bA_n^{-1}\bX^T_n \bu_n$. The consistency problem  involves finding conditions on the sequences $\bx_n$ and $u_n$ under which $\bar{\bb}_n$ converges strongly to zero. 

In \cite{drygas1976weak} and \cite{andersonTaylorAOS1976}, conditions are given for deterministic $\bX_n$ under which
$\bA_n^{-1} \rightarrow 0$ is necessary and sufficient for weak convergence of $\bar{\bb}_n$ to zero. Here it is assumed that the maximum eigenvalue of the covariance matrix of $\bu_n$ remains bounded as $n\rightarrow \infty$. In \cite{lai1979strong} it is shown that $\bA_n^{-1} \rightarrow 0$ is sufficient for strong convergence, assuming that the series $\sum_{i=1}^\infty c_i u_i$ converges $a.e.$ whenever $\sum_{i=1}^\infty c^2_i < \infty$. 

Thus, under quite general conditions, for the deterministic case weak convergence implies strong convergence.  Our purpose in this section is to extend the same general rule to the control model.  A commonly used approach is to, in essence, turn $\bA_n$ into a scalar by provisionally replacing it with a diagonal matrix $\bG_n$, the assumption then being that $\bA_n^{-1}\bG_n$  is uniformly bounded over $n$. In  \cite{anderson1979strong}, this is accomplished by assuming that the eigenvalue ratio $\lambda_{max}(\bA_n)/\lambda_{min}(\bA_n)$ is uniformly bounded over $n$, hence the sequence of matrices $\bA_n^{-1} \trace(\bA_n)$ is similarly bounded. A more general condition is used in \cite{christopeitHelmesAOS1980}, by assuming $\bA_n^{-1} \bG_n$ is uniformly bounded  over $n$, where $\bG_n = \diag(g(d^2_{n,1}), \ldots, g(d^2_{n,p}))$ for some function $g(x)$.  Here we will use this  model. Under regularity conditions to be given below,  it is verified in \cite{christopeitHelmesAOS1980} that  $\lim_n \min_t d^2_{n,t} = \infty$ is sufficient for strong consistency.  We will also add to this result by allowing  $\lim_n d^2_{n,t}  <  \infty$ for $t \leq q$ and $\lim_n d^2_{n,t}  =  \infty$ for $t > q$. Possibly, $q = 0$, in which case we will have strong consistency. Otherwise we have intermediate cases. We show that $\bb_n$ always possesses a finite limit for any $q$. Indeed, the component $\bb_n(t)$ is strongly consistent for $\beta_t$ when $t > q$.  On the other hand, if $q \geq 1$, and we define subvector $\bar{\bb}^1_n = (\bar{\bb}_n(1), \ldots, \bar{\bb}_n(q))^T$, then $P(\lim_n \bar{\bb}^1_n = \bzero) < 1$. This means that $\lim_n \min_t d^2_{n,t} = \infty$  is also necessary for strong consistency.   

Interestingly,  although \meq{linearmodel} defines a multivariate process, the main analysis  will consider univariate processes, so we make use of \mthm{thm.main.nonexpansive}.  Conditions (E1)-(E4) will define our model. Condition (E5) then suffices for strong consistency, but, as already discussed,  we will also consider intermediate cases. 
 
\begin{itemize}[noitemsep]
\item[(E1)] $u_1, u_2, \ldots$ is a sequence of $L_2$ martingale differences adapted to $\tcalf$. 
\item[(E2)] $\sup_n \var{u_n \mid \calf_{n-1}} \leq \sigma^2$  $a.e.$ for some $\sigma^2 < \infty$.   
\item[(E3)] There exists finite $n_0$ such that $\bA_{n_0}$ is nonsingular. 
\item[(E4)]  $\bA_n^{-1} \bG_n$ is uniformly bounded  over $n$, where $g(x)$ is a nondecreasing function satisfying $\int_c^\infty g^{-2}(x) dx < \infty$ for some $c > 0$. 
\item[(E5)] $\lim_n d^2_{n,t} = \infty$ for $t = 1,\ldots, p$. 
\end{itemize}
Note that (E4)  implies $\lim_{x \rightarrow\infty} g(x) = \infty$.


We may write $\bar{\bb}_n = \bA_n^{-1} \bv_n$, where 
\begin{displaymath}
\bv_n = \bX^T_n \bu_n = \bX^T_{n-1} \bu_{n-1} + u_n \bx_n = \bv_{n-1} + u_n \bx_n. 
\end{displaymath}
However, at this point we replace $\bA_n^{-1}$ with  $\bG_n^{-1}$, and define the process  
$$
\bz_n = \bG_{n}^{-1} \bv_n, \,\,\, n \geq 1. 
$$
If necessary, set $\bz_0 = \bv_0 = \bzero$. Since $\bG_{n}^{-1}$ is $\calf_{n-1}$-measurable we have  
$$
E[\bG_n^{-1} \bv_\bn \mid \calf_{n-1}] = \bG_n^{-1}  E[\bv_\bn \mid \calf_{n-1}]  = \bG_n^{-1} \bv_{n-1}.
$$
Fix $t$ and consider the component process $\bz_n(t) = \bv_n(t)/g\left(d^2_{n,t}\right)$. If $\bz_{n-1}(t) = 0$, then  $\bv_{n-1}(t)  = 0$, and hence $E[\bz_n(t) \mid \calf_{n-1}] = \bv_{n-1}(t)//g\left(d^2_{n,t}\right)  = 0$.  Therefore (A3) holds for $\bz_n(t)$.  The ratio of assumption (A1) is given by
\begin{align*}
\frac{ E[ \bz_n(t) \mid \calf_{n-1}]}{\bz_{n-1}(t)} I\{\bz_{n-1}(t) \neq 0\}
&=  \frac{g\left(d^2_{n-1,t}\right) }{ g\left(d^2_{n,t}\right) } I\{ \bz_{n-1}(t) \neq 0\}  \\
&= k_{n,t} I\{ \bz_{n-1}(t) \neq 0\}.
\end{align*}
Suppose (E4) holds.  Then $d^2_{n,t}$ is nonnegative and nondecreasing in $n$ and $g(x)$ is nondecreasing in $x$, so we may conclude that $0 \leq k_{n,t}  \leq 1$, and therefore (A1) holds for $\bz_n(t)$.  If $d^2_{n,t} \rightarrow_n \infty$, then $g(d^2_{n,t}) \rightarrow_n \infty$. It follows from \meq{eq.kronecker.01} that (A2) also holds. 

Since $\bG^{-1}_n$ is $\calf_{n-1}$-measurable,  we have residual vector 
$$
\beps_n = \bG_n^{-1} \bv_\bn - \bG_n^{-1}  E[\bv_\bn \mid \calf_{n-1}] =    \bG_n^{-1} \bx_n u_n. 
$$
This gives
\beq
E[\beps_n^T \beps_n \mid \calf_{n-1}] \leq \sigma^2 \bx_n^T \bG_n^{-2} \bx_n = \sigma^2 \sum_{t = 1}^p \frac{ x_{nt}^2 }{g^2\left(d^2_{n,t}\right)}.  \label{eq.lemma.ls.08}
\eeq
By \mlem{lem.tech.lem.01} we may conclude that $\sum_{n=1}^\infty E[\beps_n^T \beps_n \mid \calf_{n-1}] < \infty$, so that (A4) holds, following \meq{eq.durret.4.5.2}.  
Thus, by \mthm{thm.main.nonexpansive},  $\bz_n(t)$ converges to a finite limit, and if $g(d^2_{n,t}) \rightarrow_n \infty$,  $\bz_n(t)$  converges to zero $a.e.$

The convergence properties of $\bar{\bb}_n$ follow directly, by writing
\beq
 \bar{\bb}_n = \bA_n^{-1} \bv_n =   \bA_n^{-1} \bG_n \bG_n^{-1} \bv_n   =   \bA_n^{-1} \bG_n \bz_n. \label{eq.LS.46}
\eeq
Then, by (E4), $\bA_n^{-1} \bG_n$  is bounded uniformly as $n \rightarrow\infty$, so we have just proven that (E5) is sufficient for strong consistency, as reported in \cite{christopeitHelmesAOS1980}.  

If we then more precisely characterize properties of $\bA_n^{-1}$ implied by (E4) we can deduce convergence properties for intermediate cases, and show that (E5) is also necessary for strong consistency.    

\begin{theorem}\label{thm.LS.main.01}
Suppose the regression model of \meq{linearmodel} satisfies (E1)-(E4). Then the following statements hold: 
\benr
\item
The least-squares estimator $\bb_n$ possesses a finite limit $a.e.$
\item
If $\lim_n d^2_{n,t} = \infty$, then  $\bb_n(t)$ is a strongly consistent estimator of $\beta_t$. 
\item
$\bb_n$ is a strongly consistent estimator of $\bbeta$ if and only if (E5) holds. \qedhere
\een
\end{theorem}


\newpage

\appendix

\renewcommand{\thesection}{Appendix \Alph{section}}
\renewcommand{\thetheorem}{\Alph{section}.\arabic{theorem}}%
\renewcommand{\thelemma}{\Alph{section}.\arabic{lemma}}%
\renewcommand{\theequation}{\Alph{section}.\arabic{equation}}

\section{Proof of \mthm{thm.main.nonexpansive}, with Technical Lemmas}


 

\begin{lemma}\label{lem.crossing.01} 
If (A1) holds then for any $i \geq 1$ $X_i I\{D^-_{i-1} \} \leq \eps_i I\{D^-_{i-1}\}$ $a.e.$ 
\end{lemma}
\bproof We have $\left(X_i-\eps_i\right)I\{D^-_{i-1}\}  = 
    E\left[X_i \mid \calf_{i-1} \right] I\{D^-_{i-1}\} 
     \leq  0$ $a.e.$ 
by (A1), which completes the proof. 
\eproof

The following lemma allows us to bound the process $\bX$ in terms of the partial sums of the form $M_{s,t}$. Define the kernel $\lambda^{k}_t = \prod_{i=t-k+1}^t (1 + \alpha_i)$, $t \geq 1$, $k \geq 1$, and  $\lambda^{0}_t = 1$.  Under (A1), $\lambda^{max} = \sup_{t,k}  \lambda^{k}_t  < \infty$.  Then set $M^{abs}_{s,t} = \sum_{i=s}^t \absb{\eps_i}$. 

The lemma also relies on the following device. Let $\intege$ be the extended integers. For any $t \in \intege$ define $t^- = \sup \{ i \in \intege : i < t \}$.  

\begin{lemma}\label{lem.w.bound.01} 
If (A1) holds then for $j \geq 1$
\begin{align*}
W_j &\leq \lambda^{max} \left( \sup_{T_j \leq i < T_{j+1}} M^{abs}_{T_j,T_{j+1}-1}  +
\absb{U_{T_j}} \right) I\{ T_j < \infty \} \,\,\, a.e. \qedhere
\end{align*}
\end{lemma}

\bproof We have by \mlem{lem.crossing.01} and the definition of
$\eps_i$, for $i\geq 1$
\begin{eqnarray*}
X_i I\{D^-_{i-1}\} & \leq & \eps_i I\{D^-_{i-1}\} \\
X_i I\{D^0_{i-1}\} & = & (E[X_i \mid \calf_{i-1}] + \eps_i)I\{D^0_{i-1}\}
\end{eqnarray*}
which when combined give
\begin{equation}
X_i I\{(D^+_{i-1})^c \} \leq U_i + \eps_i I\{ (D^+_{i-1})^c \} = (U_i + \eps_i) I\{ (D^+_{i-1})^c \},  \label{Lemma4.1}
\end{equation}
noting that $I\{ D^0_{i-1} \} = I\{ D^0_{i-1} \} I\{ (D^+_{i-1})^c \}$. 

Fix integers $1 \leq s < t$ where $s$ is finite but, possibly, $t = \infty$. For now, assume $s < t^-$, and select $i$,  $s < i < t$. Then 
\begin{align}
X_{i} I\{D^+_{s,t^-}\}  &=  \left( \left[ X_i -  E[X_i \mid \calf_{i-1} ] \right] + E[X_i \mid \calf_{i-1}] \right)  I\{D^+_{s,t^-}\} \nonumber \\
&\leq \left( \eps_i + (1+\alpha_i) X_{i-1} \right)  I\{D^+_{s,t^-}\}, \label{eq.lemma3.01} 
\end{align}
after applying (A1). Applying  \meq{eq.lemma3.01}  iteratively gives
\begin{align}
X_{i} I\{D^+_{s,t^-}\}  &\leq  \left( \sum_{k = s+1}^i \lambda_{i}^{i-k} \eps_k + \lambda_i^{i-s} X_s \right) I\{D^+_{s,t^-}\}. \label{eq.lemma3.02} 
\end{align}
Note that \meq{eq.lemma3.02} holds also for $i = s$, using the standard summation convention, so we now permit $s=t^-$.  With \mlem{lem.crossing.01} this leads to 
\begin{align}
X_{i} I\{D^+_{s,t^-}\}I\{(D^+_{s-1})^c \} &\leq
 \left( \lambda_i^{i-s} (U_s + \eps_s)  + \sum_{k = s+1}^i \lambda_{i}^{i-k} \eps_k \right) I\{D^+_{s,t^-}\}I\{(D^+_{s-1})^c \} \nonumber \\
 & = \left( \lambda_i^{i-s}  U_s + \sum_{k = s}^i \lambda_{i}^{i-k} \eps_k \right) I\{D^+_{s,t^-}\}I\{(D^+_{s-1})^c \}   \label{Lemma4.2.1}
\end{align}
for $s \leq i < t$. It follows from \eqref{Lemma4.2.1} that 
\begin{align}
\sup_{s \leq i < t} X_i I\{D^+_{s,t^-}\}I\{(D^+_{s-1})^c \} & \leq
\lambda^{max} (U_s + M^{abs}_{s,t^-}) I\{D^+_{s,t^-}\}I\{(D^+_{s-1})^c \}.  \label{Lemma4.2.3}
\end{align}
Then on $\{ T_j < \infty\}$  we have  $\{ T_j = s, T_{j+1}=t \} \cap D^+_{s} \subset D^+_{s,t^-} \cap (D^+_{s-1})^c$. This gives
\begin{align}
\lefteqn{  \sup_{T_j \leq i < T_{j+1}} X_i I\{ T_j = s, T_{j+1}=t \}I\{ D^+_{T_j} \}  } \nonumber \\
& \leq   \lambda^{max}  (U_{T_j} + M^{abs}_{{T_j},T_{j+1}^-}) I\{T_j = s, T_{j+1}=t \}I\{ D^+_{T_j} \}. \label{Lemma4.3}
\end{align}
Summing \meq{Lemma4.3}  over $s$ and $t$ yields
\begin{align*}
W_j I\{D^+_{T_j}\} &=  \sup_{T_j \leq i < T_{j+1}} X_i I\{D^+_{T_j}\} I\{T_j < \infty \} \\
 & \leq  \lambda^{max} (\absb{U_{T_j}} + M^{abs}_{T_j,T_{j+1}^-}) I\{D^+_{T_j}\} I\{T_j < \infty \}.
\end{align*}
Then note that assumption (A1) also applies to $\{-X_i\}$, hence we may similarly conclude
\begin{displaymath}
W_j I\{D^-_{T_j}\} \leq
\lambda^{max} (\absb{U_{T_j}} + M^{abs}_{T_j,T_{j+1}^-}) I\{D^-_{T_j}\} I\{T_j < \infty \}.
\end{displaymath}
Also, by construction we have $X_i I\{D^0_{T_j}\} = 0$ for $T_j \leq i < T_{j+1}$
so that we may conclude
\begin{displaymath}
W_j \leq \lambda^{max} (\absb{U_{T_j}} + M^{abs}_{T_j,T_{j+1}^-}) I\{D^+_{T_j}\} I\{T_j < \infty \}. 
\end{displaymath}
completing the proof. \eproof

We are now in a position to prove the main results. We first consider the nonexpansive model. Given (A1), (A3) and (A4) we may claim that $\bX$ possesses a finite limit. 
It will be useful, however, to consider the case for which (A3) need not hold.  

\subsection{Proof of \mthm{thm.main.nonexpansive}}

\bproof We first consider Statement (i), that is, the case $N_T < \infty$. Fix $s \geq 1$, and define stopping time 
$$
M = \inf \{ i \geq s : X_i \leq 0 \}.
$$
Then construct the sequence
$$
X^\pr_{s+n} = (X_{s+n} - \eps_{s+n}) I\{D^+_{s+n-1}\}, \,\,\, n \geq 0,
$$
then
$$
X^*_{s+n} = X^\pr_{(s+n) \wedge M}, \,\,\, n \geq 0. 
$$
Following \mlem{lem.crossing.01} we have $X^*_{s+n} \geq 0$ for all $n \geq 0$. 
 
Then $X^*_{s+n} = 0$ for all $s+n \geq M$ and $X^*_{s+n} = X_{s+n}$ when for all $s+n < N$, noting that we may have $M = \infty$. This leads to 
\begin{align*}
E[X^*_{s+n} \mid \calf_{s+n-1} ]  &=  E[X^*_{s+n} I\{ M \geq s+n \} + X^*_{s+n} I\{ M < s+n \} \mid \calf_{s+n-1} ]  \\
&= E[X^\pr_{s+n} \mid \calf_{s+n-1} ]   I\{ M \geq s+n \}I\{D^+_{s+n-1}\} \\
&\,\,\,\,\,\,\,\,\,\, + X^\pr_M I\{ M < s+n \}I\{D^+_{s+n-1}\} \\
&= (1+\alpha_{s+n}) (X^\pr_{s+n-1} + \eps_{s+n-1}) I\{ M \geq s+n \}  \\
&\,\,\,\,\,\,\,\,\,\,  + X^\pr_{s+n-1} I\{ M < s+n \} \\ 
& \leq (1+\alpha_{s+n}) X^*_{s+n-1} +  (1+\alpha_{s+n}) \eps^+_{s+n-1}.
\end{align*}
By assumption $\sum_{n \geq 0} \alpha_{s+n} < \infty$ and  $\sum_{n \geq 0} (1+\alpha_{s+n})  < \infty$ $a.e.$, so that the almost martingale conditions of Theorem 1 of \cite{robbins1971convergence} hold.
We may then conclude that $X^*_{s+n}$ possesses a finite limit $a.e.$. This in turn implies that $X_n I\{ N_t < \infty\}$ possesses a finite limit.

We next consider Statement (ii). Since (A2) implies (A1), the conclusions of Statement (i) hold. By assumption
(A2), for $i \geq s$,
\begin{align}
(A_{s+1,i} - A_{s+1,i+1})I\{D^+_{s,\infty}\}  &=  (X_i - E[X_{i+1} \mid \calf_i]) I\{D^+_{s,\infty}\} \nonumber \\
               &\geq  (1-k_i)X_iI\{D^+_{s,\infty}\}\,\,\,  a.e.  \label{Theorem2.1}
\end{align}
If $X_iI\{D^+_{s,\infty}\}$ does not converge to 0
then there exists $c > 0$ such that for large enough $s^\prime > s$, 
$X_iI\{D^+_{s,\infty}\} \geq c$ {\it a.e.} on $D^+_{s,\infty}$ when $i \geq s^\prime$. Applying
\meq{Theorem2.1}  gives
\begin{eqnarray*}
\lim_{i\rightarrow\infty} A_{s+1,i} I\{D^+_{s,\infty}\} & = & \lim_{i\rightarrow\infty}
\left( A_{s+1,s^\prime} -  \sum_{i^\prime =s^\prime}^i (A_{s+1,i^\prime} - A_{s+1,i^\prime+1}) \right) I\{D^+_{s,\infty}\}  \\
 & \leq & \lim_{i\rightarrow\infty} \left( A_{s+1,s^\prime} -
\sum_{i^\prime = s^\prime}^i (1-k_{i^\prime})X_{i^\prime} \right) I\{D^+_{s,\infty}\} \\
    & \leq & \lim_{i\rightarrow\infty} \left( A_{s,s^\prime} -
c \sum_{i^\prime = s^\prime}^i (1-k_{i^\prime}) \right) I\{D^+_{s,\infty}\} \\
 & = & -\infty I\{D^+_{s,\infty}\} \,\, a.e.
\end{eqnarray*}
But following the proof of Statement (i) $A_{s,i}$ must possess a finite limit, hence $X_iI\{D^+_{s,\infty}\}$ must
converge to 0 $a.e.$ on  $\{N_T < \infty \} \cap \{T_{N_T} = s\} \cap D^+_s \subset D^+_{s,\infty}$.  The remaining proof is essentially the same as for Statement (i). 

We next consider Statement (iii), the case $N_T = \infty$. Since $M_{1,i}$ has a limit
{\it a.e.}, by the Cauchy criterion $M_{m,n} \rightarrow 0$ as
$m,n \rightarrow\infty$. which implies $wp1$
$$
\limsup_{j\rightarrow\infty} \left( \sup_{T_j \leq i < T_{j+1}} \absb{M_{T_j,i}} + \absb{U_{T_j}}  \right) I\{N_T = \infty\} 
= \left( \limsup_{j\rightarrow\infty} \absb{U_{j}}  \right) I\{N_T = \infty\} 
$$
since $ \{N_T = \infty \} \subset \{T_j < \infty\}$ for all $j$, and $T_j \rightarrow_j \infty$. 
Statement (iii) follows by applying \mlem{lem.w.bound.01}. Given (A3), Statement (iv) follows from a direct application of Statement (ii). 
\eproof 

\subsection{Proof of \mthm{thm.rm.01}}

\bproof Assumption (B2) ensures that the process (\ref{Algorithm}) is always defined. By (B4),  $\alpha_i\rightarrow_i 0$, so we may assume without loss of generality that $M \alpha_i \leq 1$ for all $i$. Subtracting $x^*$ from both sides of (\ref{Algorithm}) gives the equivalent process 
\begin{eqnarray}
X_i - x^* & = & X_{i-1} - x^* - \alpha_i U_i, \,\, i \geq 1, \label{eq.offset.sa.01}
\end{eqnarray}
with initial value $X_0 - x^*$. We then apply \mthm{thm.main.nonexpansive} to the process \eqref{eq.offset.sa.01}. Accordingly, by (B7), we have
\begin{displaymath}
E[X_i - x^*  \mid  \calf_{i-1}]  = X_{i-1} - x^* - \alpha_i g(X_{i-1}), \,\, i \geq 1.
\end{displaymath}
Since $g(x^*)=0$ we have
$
E[X_i - x^*  \mid  \calf_{i-1}] I\{ X_{i-1} - x^* = 0 \} = 0
$
so that (A3) holds. Then by (B3)
\beq
\frac{ E[X_i - x^*  \mid  \calf_{i-1}] }{X_{i-1}-x^*} I\{X_{i-1} - x^* \neq 0\} = \left( 1- \alpha_i \frac{g(X_{i-1})}{X_{i-1}-x^*} \right) I\{X_{i-1} - x^* \neq 0\}, \label{eq.thm.rm.01.01}
\eeq
which gives
\begin{align}
(1-M\alpha_i) I\{X_{i-1} - x^* \neq 0\} & \leq  \frac{ E[X_i - x^* \mid \calf_{i-1}] }{X_{i-1}-x^*} I\{X_{i-1} - x^* \neq 0\} \nonumber \\
& \leq  (1-m\alpha_i) I\{X_{i-1} - x^* \neq 0\}.  \label{eq.thm.rm.01.02}
\end{align}
Then note that $1 - M\alpha_i \geq 0$ and $1-m\alpha_i \leq 1$. Also, (B6) implies $\sum_i m\alpha_i = \infty$ so that (A2) holds. Then (B5) is equivalent to (A4), so that $X_i-x^*
\rightarrow_i 0$ $a.e.$  by \mthm{thm.main.nonexpansive}. 
\eproof

\section{Proofs of Theorems \ref{thm.nonuniform.contraction.01} and \ref{thm.nonuniform.contraction.02}}

\subsection{Proof of \mthm{thm.nonuniform.contraction.01}}

\bproof
Let $\tcalf$ be the filtration for the process $X_i$. For any $\tau > 0$ there exists $a.e.$ finite $N_0$ such that $\absb{\eps_i} < \tau$ for all $i \geq N_0$. Select $\tau < \delta$, then consider the offset process $X^{\delta+\tau}_{N_0+i}$, $i \geq 0$. Suppose $X^{\delta+\tau}_{N_0+i} > 0$. Since $\absb{\eps_i} < \tau$ and $X_{N_0+i} = E[ X_{N_0+i} \mid \calf_{N_0 + i-1}] + \eps_i$ we must have $X_{N_0+i} > \delta$. By (D1) we have
$
0 \leq E[ X_{N_0+i+1} \mid \calf_{N_0 + i}]  \leq X_{N_0+i} = X^{\delta+\tau}_{N_0+i} 
$.
We may then write
\begin{align*}
E[ X^{\delta+\tau}_{N_0+i+1} \mid \calf_{i}] &= E[ X_{N_0+i+1}I\{ \absb{ E[ X_{N_0+i+1} \mid \calf_{N_0 + i}] } \geq \delta + \tau \}    \mid \calf_{N_0 + i}] \\
&= E[ X_{N_0+i+1}  \mid \calf_{N_0 + i}] I\{ \absb{ E[ X_{N_0+i+1} \mid \calf_{N_0 + i}] } \geq \delta + \tau \}  \\
& \in \left[0, E[ X_{N_0+i+1}  \mid \calf_{N_0 + i}] \right].
\end{align*} 
A similar argument holds for the case $X^{\delta+\tau}_{N_0+i} < 0$, hence it follows that (A1) holds for $X_{N_0+i}^{\delta+\tau}$, $i \geq 1$. 

The proof of Statement (ii) is similar to that of Statement (i). 

To prove Statement (iii), suppose $X^{\delta+\tau}_{N_0+i-1} = 0$. First consider the case $X_{N_0+i-1} = 0$. We may write
\begin{align}
&U^{\delta+\tau}_{N_0+i}  I \{  X_{N_0+i-1} = 0 \}  \nonumber \\
&= E[X_{N_0+i} \mid \calf_{N_0 + i-1}]  I \{  X_{N_0+i-1} = 0 \} I \{ \absb{ E[X_{N_0+i} \mid \calf_{N_0 + i-1}] } \geq \delta + \tau \} \nonumber \\
&\,\,\,\,\,\,\,\,\,\, \times  I \{  X^{\delta+\tau}_{N_0+i-1} = 0 \} \nonumber \\
&= U_{N_0+i} W_{N_0+i}, \label{eq.nonuniform.contraction.01}
\end{align} 
where $W_{N_0+i} \in \{ 0, 1\}$.  Otherwise, suppose  $X_{N_0+i-1} \neq 0$ but  $\absb{E[X_{N_0+i-1} \mid \calf_{N_0+i-1} ]} \in [0,\delta+\tau)$. We must then have $\absb{X_{N_0+i-1}} \in (0, \delta+2\tau)$. 
Suppose $\absb{X_{N_0+i-1}} > \delta$. Then by (A1) $\absb{ E[X_{N_0+i} \mid \calf_{i-1}]} \leq \absb{X_{N_0+i-1}} < \delta+2\tau$. Hence, for this case
\begin{align}
\absb{U^{\delta+\tau}_{N_0+i}} I\{ \absb{X_{N_0+i-1}} > \delta \} &\leq (\delta + 2\tau).  \label{eq.nonuniform.contraction.02}
\end{align} 
Finally, suppose $\absb{X_{N_0+i-1}} \leq \delta$. By (D3)  $\absb{ E[X_{N_0+i} \mid \calf_{N_0+i-1}]} < \kappa$, hence 
\begin{align}
\absb{U^{\delta+\tau}_{N_0+i}} I\{ \absb{X_{N_0+i-1}} \in (0, \delta] \} &< \kappa. \label{eq.nonuniform.contraction.03}
\end{align} 
Then \meq{eq.thm6.(iii)} follows from Equations \eqref{eq.nonuniform.contraction.01}, \eqref{eq.nonuniform.contraction.02} and \eqref{eq.nonuniform.contraction.03}. 

We finally consider Statement (iv). Assume (A4) holds for $X_i$. Then for $X_{i}^{\delta+\tau}$ the residuals are given by 
\begin{align*}
\eps_{i}^{\delta+\tau} &= X_{i}^{\delta+\tau} - E[X_{i}^{\delta+\tau} \mid \calf_{N_0 + i-1}] \\
&=   (X_{i} - E[X_{i}])  I\{ \absb{ E[X_{N_0+i} \mid \calf_{N_0 + i-1}] } \geq \delta + \tau \} \\
&= \eps_i  I\{ \absb{ E[X_{N_0+i} \mid \calf_{N_0 + i-1}] } \geq \delta + \tau \},
\end{align*}
Thus, $\absb{ \eps_{i}^{\delta+\tau}} \leq \absb{\eps_i}$, hence (A4) extends to   $\eps_{i}^{\delta+\tau}$. 
\eproof

\subsection{Proof of \mthm{thm.nonuniform.contraction.02}}

\bproof
For the first (nonexpansive) case, we may conclude by \mthm{thm.nonuniform.contraction.01} that there exists finite $N_0$ such that $\bX^{\delta+\tau}_{N_0}$ satisfies (A1) and (A4) for all small enough $\delta, \tau > 0$. 
We then note that 
\beq
\absb{X_{N_0+i}} \leq \max\{\absb{X^{\delta+\tau}_{N_0+i}}, \delta + 2\tau\}. \label{eq.thm.nonuniform.contraction.02.01}
\eeq
 Applying  (D4) and (A3)  we conclude by  \mthm{thm.nonuniform.contraction.01} (iii) that $\limsup_i \absb{U^{\delta+\tau}_{N_0+i}} \leq \delta + 2\tau + \kappa_\delta$.   That $\bX$ possesses a finite limit follows by applying  \mthm{thm.main.nonexpansive} to $\bX^{\delta+\tau}_{N_0}$, allowing $\delta, \tau$ to approach zero, and noting \meq{eq.thm.nonuniform.contraction.02.01}.
 
For the contractive case,  note that under the given assumptions  $\absb{\bX^{\delta+\tau}_{N_0+i}}$ may be bounded by some finite constant $\delta_2$. If (D2) holds for all  $0 < \delta_1 < \delta_2 < \infty$ 
pairs then it follows from \mthm{thm.nonuniform.contraction.01} (ii) that (A2) holds for $\bX^{\delta+\tau}_{N_0}$, which therefore satisfies the hypothesis of  \mthm{thm.main.nonexpansive} (ii). The remaining argument follows that of the nonexpansive case. 
 \eproof

\section{Proofs of Theorems \ref{thm.main.multivariate.01} and  \ref{thm.rm.02}}

\subsection{Proof of \mthm{thm.main.multivariate.01}}

\bproof
Conditions (A5), (A6) and (A7) directly imply conditions (A1), (A2) and (A3) for $\bX^* = \{ \| \bX_i \| ;  i \geq 0\}$. It remains to verify that (A8) implies (A4) for $\bX^*$. We have residual process  
$$
\beps^*_i = \| \bX_i \| - E[ \| \bX_i \| \mid \calf_{i-1}], \,\,\, i \geq 1. 
$$
Then 
\begin{align}
\var{\beps^*_i \mid \calf_{i-1}} &= E[ \| \bX_i \|^2 \mid \calf_{i-1}] - E[ \| \bX_i \| \mid \calf_{i-1}]^2 \nonumber \\
&=    E\left[ \sum_{t}\bX_i(t)^2 \mid \calf_{i-1}\right] - E[ \| \bX_i \| \mid \calf_{i-1}]^2.    \label{eq.thm.main.multivariate.01}
\end{align}
However, Euclidean norm $\| \cdot \|$ is convex, so by Jensen's inequality we have 
$
E[ \| \bX_i \| \mid \calf_{i-1}]^2 \geq \sum_{t} E\left[ \bX_i(t) \mid \calf_{i-1}\right]^2,
$ 
which, when combined with  \meq{eq.thm.main.multivariate.01}, gives
\begin{align}
\var{\beps^*_i \mid \calf_{i-1}} &\leq \sum_{t}\var{\bX_i(t) \mid \calf_{i-1}} \leq \sum_{t}\var{\beps_i(t) \mid \calf_{i-1}},   \label{eq.thm.main.multivariate.02}
\end{align}
which, with (A8), implies that  $\bX^*$ satisfies (A4). The proof is completed by a direct application of  \mthm{thm.main.nonexpansive}.
\eproof

\subsection{Proof of  \mthm{thm.rm.02}}

\bproof  The main task is to verify condition (A6) for all large enough $i$.   
By construction we have
$
E[ \bX_i \mid \calf_{i-1}] = \bX_{i-1} - \alpha_i  g(\bX_{i-1}), 
$
so that, by (B3a), 
\begin{align*}
\|E[ \bX_i \mid \calf_{i-1}] \|^2 &= \| \bX_{i-1} \|^2  - 2 \alpha_i  g(\bX_{i-1}) \circ \bX_{i-1} +  \alpha_i^2  \| g(\bX_{i-1}) \|^2 \\
&\leq \| \bX_{i-1} \|^2 \left[ 1 - 2\alpha_i m + \alpha_i^2 M^2 \right].  
\end{align*}
It follows that 
\begin{displaymath}
 \frac{ \| E[\bX_i  \mid  \calf_{i-1}] \|}{ \|\bX_{i-1} \| } I\{ \|\bX_{i-1} \| \neq 0\} \leq  \left[ 1 - 2\alpha_i m + \alpha_i^2 M^2 \right]^{1/2} = k_i.  
\end{displaymath}
Since $\alpha_i \rightarrow 0$, we must have $\left[ 1 - 2\alpha_i m + \alpha_i^2 M^2 \right]^{1/2} < 1$ for all large enough $i$.   
Suppose we have function $h(u) = 1 - [1 - u]^{1/2}$. It is easily verified that $h^\pr(0) = 1/2$, and that $h(u)$ is convex for $u < 1$. Furthermore, 
there exists finite $K > 0$ and finite index $i_0$ such that for all $i \geq i_0$  we have $2\alpha_i m - \alpha_i^2 M^2 \geq K \alpha_i$. We then 
have
$
\sum_{i \geq i_0} (1 - k_i) \geq (K/2) \sum_{i \geq i_0}  \alpha_i = \infty 
$
so that (A6) holds. That (A7) and (A8) hold follow from arguments similar to those of \mthm{thm.rm.02}. The theorem is then proved following an application of  \mthm{thm.main.multivariate.01} (ii). 
\eproof

\section{Proof of \mthm{thm.LS.main.01} with Technical Lemmas}

The following lemma is a generalization of Lemma 1 of  \cite{taylor1974asymptotic}, which states that for any sequence of real numbers $x_1, x_2, \ldots$ with $x_1 \neq 0$ we must have $S_N = \sum_{n=1}^N x_n^2 / (\sum_{i=1}^n x_i^2)^2 \leq 2/x_1^2$ for all $N \geq 1$.   A recursive argument is used in \cite{taylor1974asymptotic}, but we may also interpret $S_N$ as an approximation of the integral $\int_{c}^\infty x^{-2} dx$. Doing so will allow us to expand the class of models to be studied.  

\blem\label{lem.tech.lem.01}
Let $a_1, a_2, \ldots$ be a sequence of nonnegative real numbers, not uniformly equal to zero. Assume $a_1 > 0$ (otherwise delete $a_1$ from the sequence). Let $A_n = \sum_{i=1}^n a_i$ Suppose we are given  a nondecreasing function $f:(0,\infty)\rightarrow(0,\infty)$ for which $\int_c^\infty f^{-1}(x) dx < \infty$ for any $c > 0$. Then 
\beq
S_N = \sum_{n=1}^N a_n/f(A_n) \leq a_1/A_1 + \int_{a_1}^\infty f^{-1}(x) dx \label{eq.lem.tech.lem.01.01}
\eeq
for all $N \geq 1$.
\elem

\bproof
Let $A_{max} = \lim_n A_n$. Construct a step function $h(x)$ on $x \in (0,A_{max})$ with  discontinuities at $A_n$, $n \geq 1$. Set the left limit of $h(x)$ at $A_n$ to be $f^{-1}(A_n)$. If $A_{max} < \infty$ set $h(x) = 0$ for $x \geq A_{max}$. Then 
$$
S_N = \sum_{n=1}^N a_n/f(A_n) = \int_0^{A_N} h(x) dx \leq \int_0^\infty h(x) dx = a_1/f(A_1) + \int_{a_1}^\infty h(x) dx. 
$$
\meq{eq.lem.tech.lem.01.01} follows from the fact that $0 \leq h(x) \leq f^{-1}(x)$ for $x > 0$.  
\eproof

\subsection{Proof of \mthm{thm.LS.main.01}}

For $\by = (y_1, \ldots, y_p) \in \real^p$ we have $\ell_\infty$ norm $\absb{\by}_\infty = \max_{t = 1, \ldots, p} \absb{y_t}$. Then let $\bC$ be a $p \times p$ matrix with elements $c_{ij}$. We will make use of the matrix norm $\| \bC \|_\infty = p \max_{i,j} \absb{c_{ij}}$. Note that $\| \bC \|_\infty$ is a true matrix norm, in particular,   
$\absb{\bC \by}_\infty \leq \| \bC \|_\infty  \absb{\by}_\infty$  \cite{horn2012matrix}.

\bproof
Suppose for model \eqref{linearmodel} conditions (E1)-(E4) hold. Then assume for some index $q$, $1 \leq q < p$ we have $\lim_n d^2_{n,t} <  \infty$ for $t \leq q$ and $\lim_n d^2_{n,t} = \infty$ for $t > q$.  Set $\bB_n = \bA^{-1}_n$, and let $b_{n, ij}$ be the elements of $\bB_n$. By (E4) we have $\|  \bA^{-1}_n \bG_n \|_\infty \leq \kappa$ for all $n$ for some finite $\kappa$. It follows that  
\begin{align*}
\absb{b_{n, ij}} \leq p^{-1} \kappa/g(d^2_{n,j})  \mbox{ for all } i,j.
\end{align*} 
Thus, for $j > q$ we must have $\lim_n b_{n, ij} = 0$, and, since $\bB_n$ is symmetric, for $i > q$ as well.  Then consider the matrix partitions
$$
\bA_n = \begin{bmatrix} \bA_n^{11} &  \bA_n^{12} \\  \bA_n^{21} &  \bA_n^{22} \end{bmatrix}
, \,\,\, \bB_n = \begin{bmatrix} \bB_n^{11} &  \bB_n^{12} \\  \bB_n^{21} &  \bB_n^{22} \end{bmatrix}
$$
where $\bA_n^{11}$, $\bB_n^{11}$ are $q \times q$ square matrices.
We then have
\beq
 \bA_n^{11}   \bB_n^{11} + \bA_n^{12}   \bB_n^{21} = \bI_q, \label{eq.LS.34}
\eeq
where $\bI_q$ is the $q \times q$ identity matrix. Under our given assumptions, the finite limit $\lim_n \bA_n^{11} = \bA_\infty^{11}$ exists and is positive definite.  
Let $\bC_n = \bA_n^{12}   \bB_n^{21}$ and let $c_{n,rc}$ be the elements of $\bC_n$.
Then (noting that $b_{n,ij} = b_{n,ji}$)  
\begin{align}
 \absb{c_{n,rc}} &\leq \sum_{j=q+1}^p \absb{a_{n,r j}} \absb{b_{n, jc}} \nonumber  \\
  &= \sum_{j=q+1}^p \absbb{ \sum_{i=1}^n x_{i,r} x_{n,j}} \absb{b_{n, jc}}  \nonumber   \\
  & \leq p^{-1} \kappa \left( \sum_{i=1}^\infty x^2_{i,r} \right)^{1/2} \sum_{j=q+1}^p \frac{ \left( d^2_{n,j} \right)^{1/2} }{g(d^2_{n,j})}. \label{eq.LS.36}
  \end{align}
By assumption $\sum_{i=1}^\infty x^2_{i,r} < \infty$ for $r \leq q$. Furthermore, if (E4) holds we must have $g(x) = x^{1/2}g^*(x)$, where $\lim_{x\rightarrow\infty} g^*(x) = \infty$.  By assumption $\lim_n d^2_{n,j} = \infty$ for $j >  q$, so by  \meq{eq.LS.36} it follows that $\lim_n c_{n,rc} = 0$. After applying \meq{eq.LS.34} we have $\lim_n \bB_n^{11} = \left( \bA^{11}_\infty \right)^{-1}$. 

We then have $\bb_n - \bbeta = \bA_n^{-1} \bv_n$. Define the partition  $\bv_n = \left[ \bv^1_n  \,\,\, \bv^2_n \right]$, where $\bv^1_n$ contains the first $q$ elements of $\bv_n$.  For $t \leq q$ it is easily verified that $\bv_n(t)$ is a finite variance martingale, and hence possesses limit   $\bar{\bv}(t) = \lim_n \bv_n(t)$ 
with $\var{\bar{\bv}(t)} > 0$. It follows that $P\left( \lim_n \bv^1_n = \bzero\right) < 1$. Suppose $\lim_n \bv^1_n  =  \bar{\bv}^1$.  We then have the limit 
\begin{align*}   
\lim_n \bA_n^{-1} \begin{bmatrix} \bv^1_n \\ \bzero \end{bmatrix}  &= \lim_n \begin{bmatrix} \bB^{11}_n \bv^1_n \\ \bB^{21}_n \bv^1_n \end{bmatrix}  
= \begin{bmatrix} (\bA^{11}_\infty)^{-1}  \bar{\bv}^1 \\ \bzero \end{bmatrix},  
\end{align*}
noting that the elements of $\bB^{21}_n$ vanish as $n \rightarrow\infty$. Since $\bA^{11}_\infty$ is nonsingular, $\bar{\bv}^1 \neq 0$ implies $(\bA^{11}_\infty)^{-1}  \bar{\bv}^1 \neq 0$, so that  $P\left((\bA^{11}_\infty)^{-1}  \bar{\bv}^1 = \bzero \right) < 1$. 
Next, we evaluate 
\begin{align*}   
\bW_n = \bA_n^{-1} \begin{bmatrix} \bzero \\ \bv^2_n \end{bmatrix}  &=  \bA_n^{-1} \bG_n  \bG_n^{-1}  \begin{bmatrix} \bzero \\ \bv^2_n \end{bmatrix}
\end{align*}
However, in \msec{sec.LS} it was shown that the components of index $t > q$ of $\bG_n^{-1} \bv_n$  converge to zero $a.e.$  Since $\bA_n^{-1} \bG_n$ is uniformly bounded over $n$, it follows that $\lim_n \bW_n = \bzero$ $a.e.$ We then have, $wp1$,  the finite limit
\beq
\lim_n \bb_n - \bbeta = \lim_n \bA_n^{-1} \begin{bmatrix} \bv^1_n \\ \bzero \end{bmatrix}   + \lim_n \bA_n^{-1} \begin{bmatrix} \bzero \\ \bv^2_n \bzero \end{bmatrix} 
= \begin{bmatrix} (\bA^{11}_\infty)^{-1}  \bar{\bv}^1 \\ \bzero \end{bmatrix}. \label{eq.LS.68} 
\eeq
\meq{eq.LS.68}  can clearly be extended to the case $q = p$, in which case the finite limit $\lim_n \bb_n - \bbeta = \lim_n \bA_n^{-1} \bar{\bv}$ exists $a.e.$ This is true also for condition (E5) (that is, $q = 0$), in which case $\lim_n \bb_n - \bbeta = \bzero$ $a.e.$  Hence, $\lim_n \bb_n - \bbeta$ always possesses a finite limit, and that limit is zero $wp1$ if and only if $q = 0$.  Otherwise $\bb_n(t)$ converges to $\beta_t$ $wp1$ if $\lim_n d^2_{n,t} = \infty$.   
\eproof


\end{document}